%%%%%%%%%%%%%%%%%%%%%%%%%%%%%
%% ab.tex paper on compact %%
%% generation of solvable  %%
%% subgroups of connected  %%
%% Lie groups              %%
%% Output 6 pages          %% 
%% TeXable with plain TeX  %%
\def\date{Version of       %% 
08-01-03}                  %%
\magnification=\magstep1   %% 
\hsize=15truecm            %%
\vsize=23truecm            %%
\line{\hfill\date}         %%
\parindent=1em             %%
\catcode`\@=11             %%
\newfam\msbfam             %%
\newfam\euffam             %%
\font\tenmsb=msbm10        %%
\font\sevenmsb=msbm7       %%
\font\teneuf=eufm10        %%
\font\seveneuf=eufm7       %%
\font\fiveeuf=eufm5        %%
\textfont\msbfam=\tenmsb   %%
\scriptfont\msbfam=\sevenmsb%
\textfont\euffam=\teneuf   %%
\scriptfont\euffam=%       %%
\seveneuf                  %%
\scriptscriptfont\euffam=% %%
\fiveeuf                   %%
\def\C{{\fam\msbfam C}}    %% 
\def\Q{{\fam\msbfam Q}}    %% 
\def\T{{\fam\msbfam T}}    %%
\def\R{{\fam\msbfam R}}    %%
\def\Z{{\fam\msbfam Z}}    %%
\def\g{{\fam\euffam g}}    %%
\def\h{{\fam\euffam h}}    %%
\def\b{{\fam\euffam b}}    %%
\def\r{{\fam\euffam r}}    %%
\def\z{{\fam\euffam z}}    %%
\font\eightbf=cmbx8
\font\eightrm=cmr8
\font\eightsl=cmsl8
\catcode`\@=12             %%
\def\.{{\cdot}}            %% 
\def\<{\langle}            %%
\def\>{\rangle}            %%
\def\({\big(}              %%
\def\){\big)}              %%
\def\1{{\bf1}}             %%
\def\hat{\widehat}         %%
\def\Aut{\mathop{\rm Aut}  %%
\nolimits}                 %%
\def\im{\mathop{\rm im}    %%
\nolimits}                 %%
\def\Ad{\mathop{\rm Ad}    %%
\nolimits}                 %%
\def\GL{\mathop{\rm GL}    %%
\nolimits}                 %%
\def\defi{\buildrel\rm def %%
\over=}                    %%
\def\ssk{\smallskip}       %%
\def\msk{\medskip}         %%
\def\bsk{\bigskip}         %%
\def\nin{\noindent}        %%
\def\implies{%             %%
\hbox{$\Rightarrow$}}      %%
\def\iff{%                 %%
\hbox{$\Leftrightarrow$}}  %%
          %%%%%%%%%%%%%%%%%
%%
%%%%%%MACROS FOR DIAGRAMS%%%%%%%%%%%%%%%%%%%
\def\lead{\leaders\hbox to 1.5ex{\hss${.}$\hss}\hfill}
\def\arr{\hbox to 15pt{\rightarrowfill}}
\def\larr{\hbox to 60pt{\leftarrowfill}}

\def\mapright#1{\smash{\mathop{\arr}\limits^{#1}}}

%%%%%%%%%%%%%%%%%%%%%%%%%%%%%%%%%%%%%%%%%%%%%%%%%%%%%%%%%%%
\def\sdir#1{\hbox{$\mathrel\times{\hskip -4.6pt          %% 
  {\vrule height 4.7 pt depth .5 pt}}\hskip 2pt_{#1}$}}  %%
\def\qedbox{\hbox{$\rlap{$\sqcap$}\sqcup$}}              %%
\def\qed{\nobreak\hfill\penalty250 \hbox{}               %%
\nobreak\hfill\qedbox\vskip1\baselineskip\rm}            %%
%%%%%%%%%%%%%%%%%%%%%%%%%%%%%%%%%%%%%%%%%%%%%%%%%%%%%%%%%%%
\long\def\alert#1{\smallskip\line{\hskip\parindent\vrule%
\vbox{\advance\hsize-2\parindent\hrule\smallskip\parindent.4\parindent%
\narrower\noindent#1\smallskip\hrule}\vrule\hfill}\smallskip}

\newcount\litter
\def\newlitter#1.{\advance\litter by 1 \edef#1{\number\litter}}

\newlitter\bourb.
\newlitter\hoch.
\newlitter\compbook.
\newlitter\probook.
\newlitter\morr.
\newlitter\oni.
\newlitter\rag.
\newlitter\ross.
\newlitter\yami.
\newlitter\yamii.

\vglue3cm
\centerline{\bf Solvable Subgroups of Locally Compact Groups}
\msk
\centerline{Karl Heinrich Hofmann and Karl-Hermann Neeb}
\bsk
\begingroup 
\baselineskip=5pt 

\nin{\eightbf Abstract}.\quad{\eightrm
It is shown that a closed solvable subgroup of a connected Lie group
is compactly generated. In particular, every discrete solvable subgroup
of a connected Lie group is finitely generated.  
Generalizations to locally compact groups
are discussed as far as they carry. \hfill\break
{\eightsl Mathematics Subject Classification 2000:} 22A05, 22D05, 22E15; \hfill\break
{\eightsl Key Words and Phrases:} Connected Lie group, almost connected locally
compact group, solvable subgroup, compactly generated, finitely generated.}

\endgroup

\bsk\nin
A topological group $G$ with identity component $G_0$ 
is said to be {\it almost connected} if $G/G_0$ is compact.
We shall prove the following result.
\msk
\nin
{\bf Main Theorem.} \quad
{\it A closed solvable subgroup of a locally
compact almost connected group is compactly generated.}

\msk
This result belongs to a class of ``descent'' type results that are
on record for compactly generated groups. The essay [\ross] 
 provides a good background of their history. 
It follows, in particular, that {\it a discrete solvable subgroup of an
almost connected locally compact group is finitely generated}. 
\ssk
\nin
\bf
Example S. \quad \rm The  connected
simple Lie group PSL$(2,\R)$
contains a discrete free group of infinite rank; 
such a closed subgroup  is not compactly generated.\qed 
\ssk
We remark that a  nonabelian free group is countably nilpotent (see
e.g. [\probook], Definition 10.5); that is, the descending central series
terminates at the singleton subgroup after $\omega$ steps. The Main Theorem
therefore fails for transfinitely solvable subgroups in place of solvable
ones.

\msk
The following example shows that subgroups of 
finitely generated solvable groups need not be finitely generated:
\ssk\nin\bf
Example SOL.\quad \rm
Let $\Gamma \subseteq \Q \sdir{} \Q^\times$ 
be the subgroup generated by the two elements 
$a := (0,2)$ and  $b := (1,0)$. Then
$$ \Gamma \cong \left({1\over2^\infty}\.\Z\right) \sdir{} \Z, $$
is a 2-generator metabelian group, while
the abelian subgroup ${1\over2^\infty}\.\Z\times\{0\}$ 
is not finitely generated. 
\qed
Thus, in the Main Theorem, the hypothesis ``$G/G_0$ compact'' cannot
be relaxed to ``$G/G_0$ compactly generated''.

\bsk
For {\it abelian} subgroups  the Main Theorem  
will allow us to derive a characteration theorem
for compactly generated locally compact abelian groups
as follows.

\msk
\nin
\bf Theorem. \quad \it For a locally compact  abelian  group $A$ 
the following conditions are equivalent: 

                                           {\parindent=3em

\item{\rm(1)} $A$ is compactly generated.

\item{\rm(2)} $A\cong \R^k \oplus C\oplus \Z^n$ for a unique largest compact
              subgroup $C$ and natural numbers $k$, $n$.

\item{\rm(3)} The character group $\hat A$ is a Lie group.

\item{\rm(4)} There is an almost connected  locally compact group $G$ and
              a closed subgroup $H$ such that  $A\cong H$.

}
\msk
\bf
\nin
Proof.\quad \rm
(1) \implies (2): See e.g.\ [\compbook], Theorem 7.57(ii).
  
(2) \implies (3): $\hat A \cong \hat{\R^k}\oplus \hat C\oplus \hat{\Z^n}
\cong \R^k\oplus D\oplus\T^n$ for a discrete abelian group $D$.
This is a Lie group. 

(3) \implies (2): If $\hat A$ is a Lie group, then $(\hat A)_0$ is open
and isomorphic to $\R^k\oplus \T^n$ for some $k$ and $n$; 
it is divisible, whence $\hat A
\cong (\R^k\oplus \T^n)\oplus D$ for a discrete subgroup $D$.
Hence $A=\hat{\R^k}\oplus \hat D\oplus \hat{\T^n}\cong\R^k\oplus C\oplus
\Z^n$ for the unique largest compact subgroup $C$ of $A$.

\bsk

(2) \implies (4): $A\subseteq \R^k\times C\times \R^n\cong\R^{k+n}\oplus C$,
an almost connected locally compact group.

\bsk

(4) \implies (1): Let $G$ be an almost connected 
locally compact group and
$A$ a closed abelian subgroup. Then $A$ is, in
particular, solvable. Hence the Main Theorem provides the required 
implication.\qed

By comparison with Example SOL, the situation for abelian groups
is distinctly simpler  than it is for metabelian
groups:
\msk
\nin
\bf
Corollary. \quad \rm (Morris' Theorem [\morr], [\ross])
\it A closed subgroup of a compactly generated locally compact abelian 
group is compactly generated.
\ssk
\nin\bf
Proof.\quad\rm  We proved (2)\iff(3) in the Theorem independently of the
Main  Theorem. Thus if $G$ is a locally compact compactly
generated abelian group, then $\hat G$ is an abelian Lie group.
The character group $\hat A$ of a closed subgroup $A$ of $G$,
 by duality, is a quotient of the Lie group $\hat A$ 
and thus is a Lie group. Hence $A$ is compactly generated.\qed

As we now begin a proof of the main theorem we first reduce it
 to one on connected Lie groups and its closed subgroups:
\msk 
\nin\bf
Reduction. \quad \it The {\rm Main Theorem} holds if every
closed solvable subgroup $H$ of a connected Lie group $G$ is
compactly generated.
\msk 
\nin
\bf
Proof.
\quad \rm
Indeed let $G$ be an almost connected locally compact group and 
 $N$ a compact normal subgroup such that $G/N$ is a Lie group.
The existence of $N$ is a consequence of Yamabe's Theorem saying that
each almost connected locally compact group is a pro-Lie group
([\yami,\yamii]).
Then $HN$ is a closed subgroup and $HN/N$ is a closed
solvable  subgroup $A$ of the Lie group $L=G/N$ with finitely many components.
If our claim is true for connected Lie groups $G$, then $A\cap L_0$
is compactly generated. We may assume $L=L_0A$. Then  
$A\cap L_0$ has finite index in $A$. 
Therefore $A=HN/N$ is compactly generated. Then $HN$ is compactly generated.
So $H$ is compactly generated. (See [\bourb], Chap. VII, \S 3, Lemma 3.
Also see [\ross].)\qed

This reduction allows us to concentrate on connected Lie groups $G$
and closed solvable subgroups $H$. Since any locally compact connected
group, and so in particular every connected Lie group, is compactly generated
we shall have to prove that $\pi_0(H)\defi H/H_0$ is finitely generated.
\msk 
\nin
\bf
Lemma 1.\quad \it
For a closed subgroup $H$ of a connected {\rm solvable} connected
Lie group $G$ any subgroup of  $\pi_0(H)$ is finitely generated.
\ssk
\nin
\bf
Proof. \quad\rm This is proved in [\rag], Proposition 3.8.\qed

This shows that the two generator metabelian group $\Gamma$ of 
Example SOL cannot be realized as $\pi_0(H)$ for a closed subgroup 
$H$ of a connected solvable Lie group $G$---let alone be discretely embedded
into $G$.  

\msk 
\nin
\bf
Lemma 2.\quad \it
 Let 
$$ \1 \to A \to B\mapright{q}C \to \1 $$ 
be a short exact sequence of groups. If $A$ and $C$ have the property 
that each subgroup is finitely generated, then  $B$ has this property
as well. 
\ssk\nin\bf
Proof.\quad \rm
Each subgroup $\Gamma \subseteq B$ 
is an extension of the finitely generated group $q(\Gamma)$ 
by the finitely generated group $A\cap\Gamma$, hence is finitely 
generated itself. 
\qed
\msk
\nin
\bf
Lemma 3. \quad \it
Assume that the solvable Lie group $G$ 
has the property 
that each subgroup of $\pi_0(G)$ is finitely generated. 
Let $H$ be a closed subgroup of $G$.
Then each subgroup of $\pi_0(H)$ is finitely generated. 
\ssk
\nin
\bf
Proof.\quad \rm
Let $q\colon G\to\pi_0(G)$ denote the quotient map. 
Then we have a short exact sequence 
$$ \1 \to \pi_0(H\cap G_0) \to \pi_0(H) \to q(H) \to \1. $$
As a subgroup of $\pi_0(G)$, the group $q(H)$ has the property that all 
its subgroups are finitely generated, and 
the group $\pi_0(H \cap G_0)$ has this property by Lemma 1.
Now Lemma 2 implies that each subgroup of $\pi_0(H)$ 
is finitely generated.\qed

\bf\nin
Lemma 4. \quad \it If $H$ is a closed 
solvable subgroup of $\GL_n(\C)$, then 
each subgroup of $\pi_0(H)$ is finitely generated. 
\ssk
\nin
\bf
Proof.\quad\rm
Let $S$ denote the Zariski closure of $H$. 
Then $S$ is a solvable 
linear algebraic group, so that $\pi_0(S)$ is finite
(see e.g. [\oni], Theorems 3.1.1 and 3.3.1). 
Since $H$ is a closed subgroup 
of the Lie group $S$, the assertion follows from 
Lemma 3.\qed

In order to proceed we need a further line of lemmas.
We shall call a Lie group {\it linear} if it has a faithful linear
representation. The following statement  is of 
independent interest.
\msk
\nin
\bf Proposition 5. \quad \it A connected linear Lie group has a faithful
linear representation with a closed image. 
\ssk\nin
\bf Proof. \quad \rm By [\hoch], Theorem IV.3 a connected Lie group $G$ is
 linear  if and only if it is isomorphic to
a semidiret product $B\sdir{\alpha} H$ where $B$ is a simply connected
solvable Lie group and $H$ is a linear reductive Lie group with compact
center. We set $G=B\sdir{\alpha} H$ and deduce that 
the commutator subgroup $G'$ equals $(G,B)\sdir{}(H,H)$. From
[\hoch], Theorem IV.5 it follows that $G'$ is closed in $G$. The quotient
group $G/G'$ is a direct product
$${B\over(G,B)}\times{H\over(H,H)}\cong{B\over(G,B)}\times
{Z(H)_0\over(Z(H)_0\cap (H,H))},$$
where $B/(G,B)$ is a vector group and $Z(H)_0/(Z(H)_0\cap (H,H))$
is a torus. This group has a representation mapping the vector group
$B/(G,B)$ homeomorphically on a unipotent subgroup. That is, we have
a representation $\rho\colon G\to{\rm GL}(W)$ such that
$$ \ker\rho=(G,B)H \hbox{ and } \overline{\im\rho}=\im\rho,\leqno(1)$$
the image being unipotent.

Now let $\pi\colon G\to \GL(V)$ be a faithful linear representation
and define $\zeta=\pi\oplus\rho$. We shall show that $\zeta$ has a closed
image. Suppose this is not the case.  Then
there is an $X\in\g$ such that $T\defi \overline{\zeta(\exp\R\.X)}$ is a torus
not contained in $\zeta(G)$ (see [\hoch], Proposition XVI.2.3 and Theorem 
XVI.2.4). In the Appendix we shall show that, under any representation of a
connected Lie group $G$, the commutator subgroup $G'$ has a closed image.
Thus $\zeta(G')$ is closed and $\zeta(Z(H))$ is compact since $H$ has 
a compact center. Thus $\zeta(G'Z(H))=\zeta(G')\zeta(Z(H))$ is closed
and contained in $\zeta(G)$. 
Accordingly, $X$ cannot be contained in $\g'+\z(\h)=[\g,\b]+\h$.
Thus by (1), $\exp \R\.X$ fails to be in $\ker\rho$. It follows that 
$\rho \circ\exp$ maps $\R\.X$ homeomorphically onto a unipotent one-parameter
group. Then $\zeta\circ\exp$ maps $\R\.X$ homeomorphically as well, and
that contradicts the fact that $T$ is a torus. This contradiction
proves the proposition.\qed  

We now complete the proof of the Main Theorem by proving
the last lemma:

\msk
\nin
\bf
Lemma 6.\quad \it Let $G$ be a connected Lie group and $H$ a closed
solvable subgroup. Then $H$ is compactly generated.
\ssk
\nin
\bf Proof. \quad \rm Let 
$Z=Z(G)$ be the center of $G$. Then $A\defi
\overline{ZH}$ is a closed solvable subgroup of $G$ containing $H$.
By Lemma 3 for $H$ to be compactly generated it will
suffice to show that all subgroups of  $\pi_0(A)=A/A_0$ are
finitely generated. Let $A_1$ be a subgroup of $A$ containing $A_0$.
Then $A_1$ is open in $A$, and so $A_1Z$ is open and thus closed
in $A$. Therefore 
$$A_1/(A_1\cap(A_0Z))\cong A_1Z/A_0Z.\leqno(1)$$
By the modular law,
$$A_1\cap(A_0Z)=A_0(A_1\cap Z).\leqno(2)$$
We have the following isomorphism of discrete groups
$$ A_0(A_1\cap Z)/A_0\cong (A_1\cap Z)/(A_0\cap(A_1\cap Z))
= (A_1\cap Z)/(A_0\cap Z).\leqno (3)$$
Taking (1), (2) and (3) together we recognize the
following exact sequence
$$\1\to{A_1\cap Z\over A_0\cap Z}\to{A_1\over A_0}
    \to{A_1Z\over A_0Z}\to\1.\leqno(4)$$
In order to show that $A_1/A_0$ is finitely generated it
therefore suffices that 

{\parindent=3em

\item{(a)} $(A_1\cap Z)/(A_0\cap Z)$ is finitely generated,
\item{(b)} $(A_1Z)/(A_0Z)$ is finitely generated.

}
\ssk\nin
Ad (a): The center $Z$ of the connected Lie group $G$ is 
compactly generated. (Indeed the 
fundamental group $\pi_1(G/Z)$ is finitely
generated abelian and $\pi_0(Z)=Z/Z_0$ is the kernel of the covering
morphism $G/Z_0\to G/Z$ and is therefore finitely generated as a 
quotient of $\pi_1(G/Z)$. Thus $Z$ is compactly generated.)
%By Morris' Theorem 
Since $A_1$ is open in $A$, the group $A_1\cap Z$ is open in $Z$ and
thus compactly generated, and so (a) follows.
\ssk\nin
Ad (b): The adjoint representation  
$\Ad\colon G\to\Aut\g\subseteq \GL(\g)$
induces a faithful linear representation of $G/Z$. 
Then by Lemma 4 and
Proposition 5, $A_1Z/Z$, a closed solvable subgroup of $G/Z$,
is compactly generated. Then the discrete factor group
$A_1Z/A_0Z\cong(A_1Z/Z)/(A_0Z/Z)$ is finitely generated.
Thus (b) is proved as well and this completes the proof 
of Lemma 6 and thereby the proof of the Main Theorem.\qed

\centerline{\bf Appendix}
\msk\nin
In the proof of Proposition 5 we used the following 
\ssk
\nin
\bf Theorem A. \quad \it For any finite dimensional representation of a
connected Lie group $G$, the image of the commutator subgroup is closed.
\bsk
\nin
\bf Proof. \quad \rm It is no loss of generality to assume that $G$ is
simply connected. Then we have Levi decomposition $G =R\sdir{\alpha}S$
and $G'=(G,R)\sdir{} S$. Let $\pi\colon G\to {\rm GL}(V)$ be a finite 
dimensional representation and let
$$V_0=\{0\}\subseteq V_1\subseteq\cdots\subseteq V_n=V$$
be a maximal flag of $G$-submodules of $V$ such that all quotient modules
$V_{j+1}/V_j$ are simple. Since $\pi|S$ is a semisimple representation, 
we may choose $S$-invariant decompositions $V_j=V_{j-1}\oplus W_j$.
Then
$$\pi(G)\subseteq G_F\defi \{g\in{\rm GL}(V): (\forall j)gV_j=V_j\},$$
and we have a semidirect decomposition $G_F=U_F\sdir{}L_F$, where
$$U_F=\{g\in {\rm GL}(V):(\forall j)(g-1)(V_j)=V_{j-1}$$ and
$L_F=\prod_j{\rm GL}(W_j)$. Note also that $\pi(S)\subseteq L_F$.
Furthermore, Theorem I.5.3.1 of [\bourb] implies that the ideal $[\g,\r]$
acts trivially on each simple $\g$-module and so $\pi((G,R))\subseteq
U_F$. Hence $\pi((G,R))$ is a unipotent analytic group and is therefore
closed. Moreover, $\pi(S)$ is closed (see [2], Chapter XVI) 
and this shows that $\pi(G')\cong\pi((G,R))\sdir{}\pi(S)$ is closed.\qed 

The proof of Theorem A can be derived from the theory of algebraic groups,
since the commutator algebra of a linear Lie algebra is the Lie algebra
of an algebraic group [\oni]. We gave a more direct proof inspired by 
the discussion of linear Lie groups in [\hoch].

\bsk
\centerline{\bf References}
\msk
{\parindent2em

\item{[\bourb]} Bourbaki, N., Groupes et alg\`ebres de Lie, Chap.~I-III, 
reprinted by Springer-Verlag, Berlin etc., 1989. 
\ssk
\item{[\hoch]} Hochschild, G., The Structure of Lie Groups, Holden Day, 
San Francisco, 1965.
\ssk
\item{[\compbook]} Hofmann, K. H. and S. A. Morris, The Structure of 
Compact Groups, W. DeGruyter, Berlin 1998 and 2006.
\ssk
\item{[\probook]} ---, The Lie Theory of Connected Pro-Lie Groups, 
European Mathematical Society Publishing House, Z\"urich, 2007.
\ssk
\item{[\morr]}  Morris, S. A., Locally compact abelian groups 
and the variety of
topological groups generated by the reals, Proc. Amer. Math. Soc. 
{\bf34} (1972), 290--292.
\ssk
\item{[\oni]} Onishchik, A. L., and E. B. Vinberg, Lie Groups and Algebraic
Groups, Springer-Verlag, Berlin etc., 1990.
\ssk
\item{[\rag]} Raghunathan, M. S., 
``Discrete Subgroups of Lie Groups,'' Ergebnisse der Math.\ {\bf 68}, 
Springer, Berlin etc., 1972. 
\ssk
\item{[\ross]} Ross, K., Closed subgroups of compactly generated LCA 
group are compactly generated,
{\tt http://www.uoregon.edu/~ross1/subgroupsofCGLCA6.pdf}.
\ssk

\item{[\yami]} Yamabe, H.,
         {\it On the Conjecture of Iwasawa and Gleason},
             Ann. of Math. {\bf58} (1953), 48--54.
\ssk
\item{[\yamii]} ---,
         {\it Generalization of a theorem of Gleason}, 
             Ann. of Math. {\bf58} (1953), 351--365.

}

\bsk 
\baselineskip=5pt 
\parindent0pt
\eightrm
\obeylines
Karl Heinrich Hofmann
Karl-Hermann Neeb
Fachbereich Mathematik
Technische Universit\"at Darmstadt
Schlossgartenstrasse 7
64289 Darmstadt
hofmann@mathematik.tu-darmstadt.de
neeb@mathematik.tu-darmstadt.de
 
\bye